\documentclass[10pt]{article}

\def\CMP{Commun.~Math.~Phys.~}
\def \F {\mathcal{F}}
\def \D {\hbox{d}}

\def \bfv {{\bf v}}
\def\PVI    {{\rm P6}}
\def\irrec{k}

\begin{document}

\title{Symmetry reductions of a particular set of 
equations of associativity
in twodimensional topological field theory\footnote
{Journal of Physics A, to appear. Corresponding author RC. Preprint S2004/045.}
}

\author{Robert Conte\dag\ and Maria Luz Gandarias\ddag
{}\\
\\ \dag Service de physique de l'\'etat condens\'e (URA 2464), CEA--Saclay
\\ F--91191 Gif-sur-Yvette Cedex, France
\\ E-mail:  Conte@drecam.saclay.cea.fr
{}\\
\\ \ddag Departamento de Matematicas
\\ Universidad de C\'adiz
\\ Casa postale 40
\\ E--11510 Puerto Real, C\'adiz, Spain
\\ E-mail:  MariaLuz.Gandarias@uca.es
}

\maketitle

\hfill 

{\vglue -10.0 truemm}
{\vskip -10.0 truemm}

\begin{abstract}
The WDVV equations of associativity arising in twodimensional
topological field theory can be represented, in the simplest nontrivial case, 
by a single third order equation of the Monge-Amp\`ere type. 
By investigating its Lie point symmetries, 
we reduce it to various nonlinear ordinary differential equations,
and we obtain several new explicit solutions.
\end{abstract}




\noindent \textit{PACS 1995}
  02.20.Qs   
  11.10.Lm   

\baselineskip=12truept 


\tableofcontents

\vfill \eject

\section{Introduction}

As introduced by
Witten, Dijkgraaf, H.~Verlinde and E.~Verlinde \cite{W1990,DVV1991},
the equations of associativity involve the following unknowns:
a function $\F(t^1,\dots,t^n) \equiv \F(t)$,
integer numbers $q_\alpha$ and $r_\alpha$, $\alpha=1,\dots,n$,
another integer $d$,
a constant symmetric nondegenerate matrix $(\eta^{\alpha \beta})$,
other constants $A_{\alpha \beta}, B_\alpha,C$.
These unknowns must obey three main sets of equations 
\cite{Cargese1996Dubrovin}.
\begin{enumerate}
\item
the equations of associativity properly said
(with summation over the repeated indices)
\begin{eqnarray}
& &
\partial_\alpha \partial_\beta  \partial_\lambda \F(t) \eta^{\lambda \mu}
\partial_\mu    \partial_\gamma \partial_\delta  \F(t)
=
\partial_\delta \partial_\beta  \partial_\lambda \F(t) \eta^{\lambda \mu}
\partial_\mu    \partial_\gamma \partial_\alpha  \F(t)
{\hskip 23.0 truemm}  \hbox{(WDVV1)}
\nonumber
\end{eqnarray}

\item
a condition singling out one variable, say, $t^1$,
\begin{eqnarray}
& &
\partial_\alpha \partial_\beta \partial_1 \F(t) = \eta^{\alpha\beta},
{\hskip 79.0 truemm}  \hbox{(WDVV2)}
\nonumber
\end{eqnarray}
in which the matrix $(\eta_{\alpha\beta})$ is the inverse of
$(\eta^{\alpha\beta})$,

\item
a condition of quasi-homogeneity,
\begin{eqnarray}
& &
\sum_{\alpha=1}^n
\left\lbrack (1-q_\alpha)t^\alpha+r_\alpha \right\rbrack \partial_\alpha \F(t)
=
(3-d) \F(t) + \frac{1}{2} A_{\alpha\beta} t^\alpha t^\beta
 + B_\alpha t^\alpha +C.
{\hskip 17.0 truemm}  \hbox{(WDVV3)}
\nonumber
\end{eqnarray}

\end{enumerate}

In the simplest nontrivial case $n=3$,
there essentially exist two different choices of coordinates \cite{DubLNM},
depending on $\eta^{11}$ being zero or nonzero,
each choice resulting in a representation of the generating function
$\F$ in terms of the solution of 
a single third order partial differential equation (PDE)
of the Monge-Amp\`ere type,
which is either
\cite{DubLNM}
\cite[Eq.~(9)]{F2004}
\begin{enumerate}
\item
\begin{eqnarray}
& &
\eta^{11} \not=0\ :\
\F=\frac{1}{6} \left(t^1\right)^3 + t^1 t^2 t^3 + f(t^2,t^3),
\\
& &
f_{xxx} f_{yyy} -f_{xxy} f_{yyx} -1 =0,\
x=t^2,\
y=t^3,
\label{eqPDEFerapontov}
\end{eqnarray}
or
\cite[page 304]{Cargese1996Dubrovin} \cite[Eq.~(22)]{F2004}
\item
\begin{eqnarray}
& &
\eta^{11} =0\ :\
\F=\frac{1}{2} \left(t^1\right)^2 t^3 + \frac{1}{2} t^1 (t^2)^2+ F(t^2,t^3),
\\
& & \left(F_{tyy}\right)^2 - F_{ttt} -F_{tty} F_{yyy}=0,\
y=t^2,\
t=t^3.
\label{eqPDEDubrovin}
\end{eqnarray}

\end{enumerate}

There exists a Legendre transformation \cite{DubLNM}
which exchanges these two solutions $\F$
(this transformation exchanges the coordinates
$t^3$ of (\ref{eqPDEFerapontov}) and $t^2$ of (\ref{eqPDEDubrovin}),
whose common value is here denoted $y$),
and its action on the functions of two variables $f(x,y)$ and $F(y,t)$
is the hodograph transformation \cite{FM1996}
\cite[Eq.~(23)]{F2004}
\begin{eqnarray}
& &
t=f_{xx},\
F_{yyy}= \frac{f_{xxy}^2}{f_{xxx}},\
F_{tyy}=-\frac{f_{xxy}}{f_{xxx}},\
F_{tty}= \frac{1}{f_{xxx}},\
F_{ttt}= \frac{f_{xyy}}{f_{xxx}},\
\label{eqHodographFrom_f_to_F}
\end{eqnarray}
whose inverse is
\begin{eqnarray}
& &
f_{xx}=t,\
f_{xy}=-F_{yy},\
f_{yy}= F_{tt},\
x=F_{ty}.
\label{eqHodographFrom_F_to_f}
\end{eqnarray}
A nice way to obtain this hodograph transformation
(\ref{eqHodographFrom_f_to_F})--(\ref{eqHodographFrom_F_to_f})
is to rewrite both PDEs \cite{FM1996}
as integrable systems of the so-called hydrodynamic type,
allowing them to be mapped by a chain of standard transformations
to integrable three-wave systems.

Both PDEs admit a Lax pair \cite{DubLNM},
e.g.~for the PDE (\ref{eqPDEFerapontov})
\cite[Eq.(10)]{F2004}
\begin{eqnarray}
& &
\psi_x=\lambda \pmatrix {0 & 1 & 0 \cr
                         0 & f_{xxy} & f_{xxx} \cr
                         1 & f_{xyy} & f_{xxy} \cr} \psi,\
\psi_y=\lambda \pmatrix {0 & 0 & 1 \cr
                         1 & f_{xyy} & f_{xxy} \cr
                         0 & f_{yyy} & f_{xyy} \cr} \psi,
\label{eqLaxFerapontov}
\end{eqnarray}
in which $\lambda$ is a nonzero spectral parameter.

The purpose of this paper is to obtain new explicit solutions of
either the PDE (\ref{eqPDEFerapontov})
or the PDE (\ref{eqPDEDubrovin}),
and therefore of the equations of associativity
in the simplest nontrivial case.
Any such solution $f$ is only defined up to an
arbitrary additive second degree polynomial.
However,
the equation (\ref{eqPDEDubrovin})
possesses a rather complicated structure of singularities,
making uneasy the search for explicit solutions,
while the equation (\ref{eqPDEFerapontov}) has a simpler such structure,
so we will mainly consider this latter equation.
In particular,
the invariance of this PDE under permutation of $x$ and $y$ has no simple
equivalent for the PDE (\ref{eqPDEDubrovin}).

To achieve this search for solutions,
we perform a systematic investigation,
\textit{via} the Lie point symmetries method,
of the reductions of the PDE (\ref{eqPDEFerapontov}) to
ordinary differential equations (ODEs),
which \textit{a priori} can be integrated since they inherit
the integrability properties of the equations of associativity.
In addition to the reductions or particular solutions of either
(\ref{eqPDEDubrovin}) or (\ref{eqPDEFerapontov})
which have already been found \cite{DubLNM,F2004},
we obtain several new results. 

The paper is organized as follows.
In section (\ref{sectionClassicalSymmetries}),
we apply the classical Lie method \cite{OvsiannikovBook,OlverBook},
derive the Lie algebra,
compute the commutator table and the adjoint table \cite{OlverBook},
which then allow us to derive the optimal system of generators.
In section \ref{sectionClassicalReductions},
we perform all the associated classical reductions.
The last section (\ref{sectionSummary}) summarizes the solutions.

\section{Classical Lie symmetries}
\label{sectionClassicalSymmetries}

 In order to apply the classical Lie method to the Ferapontov equation
(\ref{eqPDEFerapontov}),
we consider the one-parameter Lie group of infinitesimal transformations in
$(x,y,f)$
\begin{equation}
\label{tran}
\begin{tabular}{l}
$x^*=x+\varepsilon  \xi(x,y,f)+{\mathcal O}(\varepsilon^2)$,\\[5pt]
$y^*=y+\varepsilon \eta(x,y,f)+{\mathcal O}(\varepsilon^2)$,\\[5pt]
$f^*=f+\varepsilon \phi(x,y,f)+{\mathcal O}(\varepsilon^2)$,\\[5pt]
\end{tabular}
\end{equation}
where $\varepsilon$ is the group parameter.
 The associated Lie algebra of infinitesimal symmetries is the set
of vector fields of the form

\begin{equation}
\label{vect}
\bfv=\xi \partial_x +\eta \partial_y+\phi \partial_f.
\end{equation}

One then requires that this
transformation leaves invariant the set of solutions of the
equation (\ref{eqPDEFerapontov}).
This yields an overdetermined, linear
system of equations for the infinitesimals
$\xi(x,y,f)$, $\eta(x,y,f)$ and $\phi(x,y,f)$.
 Having determined the infinitesimals, the
symmetry variables are found by solving the invariant surface condition
\begin{equation}
\label{sur}
\Phi \equiv \xi\frac{\partial f}{\partial x}
         + \eta\frac{\partial f}{\partial y}-\phi=0.
\end{equation}

Applying the classical method to the equation (\ref{eqPDEFerapontov})
leads to a ten-parameter Lie group.
Associated with this Lie group we have a Lie algebra which can be
represented by the following generators:
$$\begin{array}{llllll}
\bfv_1=    \partial_x,\hspace{0.3mm}&
\bfv_2=    \partial_y,\hspace{0.3mm}&
\bfv_3=x   \partial_x +\frac{3}{2} f \partial_f, \hspace{0.3mm}&
\bfv_4=y   \partial_y +\frac{3}{2} f \partial_f, \hspace{0.3mm}&
\bfv_5=x y \partial_f,\hspace{0.3mm}&
\cr \\[10pt]
\bfv_6=x^2 \partial_f,\hspace{0.3mm}&
\bfv_7=y^2 \partial_f,\hspace{0.3mm}&
\bfv_8=x   \partial_f,\hspace{0.3mm}&
\bfv_9=y   \partial_f,\hspace{0.3mm}&
\bfv_{10}= \partial_f.&
\cr
\end{array}
$$

\subsection{Optimal system}

In order to construct the optimal system, following Olver\cite{OlverBook},
we first construct
the commutator table (Table~\ref{TableCommutator}) and the adjoint table
(Table~\ref{TableAdjoint}) which shows the separate adjoint actions of
each element in $\bfv_i$, $i=1\dots10$, as it acts on all other elements.
This construction is done easily by summing the Lie series.

The corresponding generators of the optimal system of subalgebras are
\begin{eqnarray}
            &&                     \bfv_3,
\nonumber\\ &&                              \bfv_4,\
\nonumber\\ &&                  -a \bfv_3+b \bfv_4,
\nonumber\\ &&                     \bfv_3-  \bfv_4+a \bfv_5                               +b \bfv_{10},
\nonumber\\ &&                   3 \bfv_3+  \bfv_4                  +a \bfv_7,
\nonumber\\ &&                     \bfv_3+3 \bfv_4         +a \bfv_6,
\nonumber\\ &&                  -3 \bfv_3+  \bfv_4                           +a \bfv_8,
\nonumber\\ &&                     \bfv_3-3 \bfv_4                                    +a \bfv_9,
\nonumber\\ &&          a \bfv_2+b \bfv_3,
\nonumber\\ && a \bfv_1                  +b \bfv_4,
\nonumber\\ && a \bfv_1+b \bfv_2                  +c \bfv_5+d \bfv_6+e \bfv_7,\
\label{eqOptimalSystemGenerators}
\end{eqnarray}
where $a,b,c,d,e$ are arbitrary real nonzero constants.

\vfill \eject
\tabcolsep=0.5truemm

\vspace{20pt}
\begin{table}[ht]
\caption{Commutator table for the Lie algebra $\bfv_i$.
\label{TableCommutator}}
\begin{center}
\footnotesize
\begin{tabular}{|l|llllllllll|}
\hline & $\bfv_1$ & $\bfv_2$ & $\bfv_3$ & $\bfv_4$ & $\bfv_5$ & $\bfv_6$
       & $\bfv_7$ & $\bfv_8$ & $\bfv_9$ & $\bfv_{10}$
\\[10pt]\hline
$\bfv_1$ & $0$ & $0$ & $\bfv_1$ & $0$ & $\bfv_9$ & $2 \bfv_8$ & $0$ & $\bfv_{10}$  & $0$ & $0$
\\[10pt]
 $\bfv_2$ & $0$ & $0$ & $0$ & $\bfv_2$ & $\bfv_8$ &  $0$ &  $2\bfv_9$ & $0$ & ${\bfv}_{10}$ & $0$
 \\[10pt]
$\bfv_3$ &  $-\bfv_1$ & $0$ &  $0$ &  $0$ &
$-\frac{1}{2}\bfv_5$ & $\frac{1}{2}\bfv_6$
& $-\frac{3}{2}\bfv_7$  & $-\frac{1}{2}\bfv_8$  & $-\frac{3}{2}\bfv_9$  & $-\frac{3}{2}\bfv_{10}$\\[10pt]
$\bfv_4$ & $0$ &  $-\bfv_2$ &  $0$& $0$ & $-\frac{1}{2}\bfv_5$ & $-\frac{3}{2}\bfv_6$ & $\frac{1}{2}\bfv_7$ &
$-\frac{3}{2}\bfv_8$ &
$-\frac{1}{2}\bfv_9$ &  $-\frac{3}{2}\bfv_{10}$   \\[5pt]
$\bfv_5$ & $-\bfv_9$ & $-\bfv_8$ & $\frac{1}{2}\bfv_5$
&
$\frac{1}{2}\bfv_5$  &  $0$ &  $0$ &  $0$& $0$ &  $0$& $0$   \\[5pt]
 $\bfv_6$ & $-2\bfv_8$ & $0$  & $-\frac{1}{2}\bfv_6$ &
$\frac{3}{2}\bfv_6$  &  $0$ &  $0$ &  $0$& $0$ &  $0$& $0$   \\[5pt]
$\bfv_7$ & $0$ & $-2\bfv_{9}$  &  $\frac{3}{2}\bfv_7$ &
$-\frac{1}{2}\bfv_7$  &  $0$ &  $0$ &  $0$ & $0$ &  $0$  &  $0$   \\[5pt]
$\bfv_8$ &  $-\bfv_{10}$ & $0$ & $\frac{1}{2}\bfv_8$ &
$\frac{3}{2}\bfv_8$  &  $0$ &  $0$ &  $0$& $0$ &  $0$& $0$   \\[5pt]
$\bfv_9$ & $0$ & $-\bfv_{10}$  & $\frac{3}{2}\bfv_9$ &
$\frac{1}{2}\bfv_9$  &  $0$ &  $0$ &  $0$& $0$ &  $0$& $0$   \\[5pt]
$\bfv_{10}$ & $0$ & $0$ &  $\frac{3}{2}\bfv_{10}$ &
$\frac{3}{2}\bfv_{10}$  &  $0$ &  $0$ &  $0$& $0$ &  $0$& $0$   \\[5pt]
\hline
\end{tabular}
\end{center}
\end{table}

\begin{table}[ht] 
\caption{Adjoint table for the Lie algebra $\bfv_i$.
\label{TableAdjoint}
}
\begin{center}
\footnotesize
\begin{tabular}{|l|llllllllll|}
\hline $Ad$ & $\bfv_1$ & $\bfv_2$ & $\bfv_3$ & $\bfv_4$ &$\bfv_5$ & $\bfv_6$
& $\bfv_7$ & $\bfv_8$ & $\bfv_9$ & $\bfv_{10}$\\[10pt]\hline
$\bfv_1$ & $\bfv_1$  & $\bfv_2$  & $\bfv_3-\varepsilon
\bfv_1$ & $\bfv_4$ & $\bfv_5-  \varepsilon \bfv_9$ &
$\bfv_6-2\varepsilon \bfv_8+\varepsilon^2 \bfv_{10}$
& $\bfv_7$ & $\bfv_8-\varepsilon \bfv_{10}$ &  $\bfv_9$ &$\bfv_{10}$ \\[10pt]
 $\bfv_2$ & $\bfv_1$ & $\bfv_2$ & $\bfv_3$ & $\bfv_4-\varepsilon \bfv_{2}$
& $\bfv_5-\varepsilon \bfv_{8}$ & $\bfv_6$ & $\bfv_7-2\varepsilon \bfv_{9}+\varepsilon^2 \bfv_{10}$ & $\bfv_8$ & $\bfv_9-\varepsilon \bfv_{10}$ & $\bfv_{10}$ \\[10pt]
$\bfv_3$ &  $e^{\varepsilon}\bfv_1$ & $\bfv_2$ & $\bfv_3$
& $\bfv_4$  & $e^{\frac{\varepsilon}{2}}\bfv_5$ &
$e^{-\frac{\varepsilon}{2}}\bfv_6$ & $e^{\frac{3\varepsilon}{2}}\bfv_7$ & $e^{\frac{\varepsilon}{2}}\bfv_8$ &
$e^{\frac{3\varepsilon}{2}}\bfv_9$ &$e^{\frac{3\varepsilon}{2}}\bfv_{10}$
\\[10pt]
$\bfv_4$ & $\bfv_1 $ &  $e^{\varepsilon}\bfv_2$ & $\bfv_3$ & $\bfv_4$  & $e^{\frac{ \varepsilon}{2}}\bfv_5$ &
$e^{\frac{3 \varepsilon}{2}}\bfv_6$ & $e^{-\frac{
\varepsilon}{2}}\bfv_7$ & $e^{\frac{3 \varepsilon}{2}}\bfv_8$ &
$e^{\frac{ \varepsilon}{2}}\bfv_9$&
$e^{\frac{3 \varepsilon}{2}}\bfv_{10}$    \\[5pt]
$\bfv_5$ & $\bfv_1+\varepsilon \bfv_9$ & $\bfv_2+\varepsilon
\bfv_8$ & $\bfv_3-\frac{1}{2}\varepsilon \bfv_5$ & $\bfv_4-\frac{1}{2}\varepsilon \bfv_5$ & $\bfv_5$ &
$\bfv_6$  & $\bfv_7$  & $\bfv_8$ & $\bfv_9$ & $\bfv_{10}$   \\[5pt]
$\bfv_6$ & $\bfv_1+2\varepsilon \bfv_8$ & $\bfv_2$ &
$\bfv_3+\frac{1}{2}\varepsilon \bfv_6$ & $\bfv_4-\frac{3}{2}\varepsilon \bfv_6$ & $\bfv_5$ & $\bfv_6$ &
$\bfv_7$  & $\bfv_8$ & $\bfv_9$ & $\bfv_{10}$
\\[5pt]
$\bfv_7$ & $\bfv_1$ & $\bfv_2+2\varepsilon \bfv_9$  &
$\bfv_3-\frac{3}{2}\varepsilon \bfv_7$ & $\bfv_4+\frac{1}{2}\varepsilon \bfv_7$ & $\bfv_5$ & $\bfv_6$ &
$\bfv_7$  & $\bfv_8$ & $\bfv_9$ & $\bfv_{10}$
\\[5pt]
$\bfv_8$ & $\bfv_1+\varepsilon \bfv_{10}$ & $\bfv_2$ &
$\bfv_3-\frac{1}{2}\varepsilon \bfv_8$ & $\bfv_4-\frac{3}{2}\varepsilon \bfv_8$  & $\bfv_5$ & $\bfv_6$ &
$\bfv_7$  & $\bfv_8$ & $\bfv_9$ & $\bfv_{10}$
\\[5pt]
$\bfv_9$ & $\bfv_1$ & $\bfv_2+\varepsilon \bfv_{10}$  &
$\bfv_3-\frac{3}{2}\varepsilon \bfv_9$ & $\bfv_4-\frac{1}{2}\varepsilon \bfv_9$
& $\bfv_5$ & $\bfv_6$ & $\bfv_7$  & $\bfv_8$ & $\bfv_9$ & $\bfv_{10}$
\\[5pt]
$\bfv_{10}$ & $\bfv_1$ & $\bfv_2$ & $\bfv_3-\frac{3}{2}\varepsilon \bfv_{10}$ & $\bfv_4-\frac{3}{2}\varepsilon \bfv_{10}$ & $\bfv_5$ & $\bfv_6$
& $\bfv_7$  & $\bfv_8$ & $\bfv_9$ & $\bfv_{10}$ \\[5pt]\hline
\end{tabular}
\end{center}
\end{table} 

\vfill \eject

\section{Classical reductions}
\label{sectionClassicalReductions}

Each generator of the optimal system defines a reduction of the equation
(\ref{eqPDEFerapontov}) to an ODE.
Because of the invariance of (\ref{eqPDEFerapontov})
under permutation of $x$ and $y$,
these ten generators only define seven different reductions to an ODE,
which we now consider.

Although these reductions are probably integrable in some sense,
performing their explicit integration is a difficult task.
Moreover, since the Lax pair (\ref{eqLaxFerapontov}) is not isospectral,
its reductions, which are also Lax pairs for the reduced ODEs,
cannot generate any first integral,
so the Lax pair is unfortunately of no use for integrating the reduced ODEs.

{}From the scaling invariance of the two considered PDEs,
an obvious solution is
\begin{eqnarray}
& &
f=2 i \frac{\sqrt{2}}{3} (x y)^{3/2},\
F=\frac{y^4}{8 t},\
2 t^2 x + y^3=0,\
i^2=-1.
\label{eqScalingSolution}
\end{eqnarray}

\subsection{Reduction with the generator $ \bfv_3$ or $ \bfv_4$}

The generators $ \bfv_3$ and $\bfv_4$ define a reduction to the same
autonomous linear ODE \cite[Eq.~(30) p.~46]{F2004},
\begin{eqnarray}
& & \left\lbrace \matrix{
\displaystyle{
 z=y,\  f=\left[x^3 \Phi(z)\right]^{1/2},
 \hbox{ or }
 z=x,\  f=\left[y^3 \Phi(z)\right]^{1/2},
 }
\hfill \cr \displaystyle{
\Phi''' + 16/3=0.
} \hfill \cr} \right.
\label{eqoptim1or2}
\end{eqnarray}
This contains the scaling solution (\ref{eqScalingSolution}).

\subsection{Reduction with the generator $- a \bfv_3 + b \bfv_4$}

With the notation $s=a+b,p=ab$, a symmetric definition of this reduction is,
\begin{eqnarray}
& & \left\lbrace
\begin{array}{ll}
\displaystyle{
 z=x^b y^a,\ f=(x y)^{3/2} \varphi(z),
}
\\
\displaystyle{
\left[
-16 p^2 s z^5 \varphi''
-8 p (4 p +2 p s + s^2) z^4 \varphi'
- 3 s^3 z^3 \varphi
\right] \varphi'''
} \\ \displaystyle{ \phantom{1234}
+ 8 p (2 p - 6 p s - 3 s^2) {\varphi''}^2
- (64 p^2 +72 p s +64 p^2 s +72 p s^2 +9 s^3) z^3 \varphi' \varphi''
} \\ \displaystyle{ \phantom{1234}
- 9 (2 +s) s^2 z^2 \varphi  \varphi''
- (40 p +16 p^2 +72 p s +16 p^2 s +18 s^2 +32 p s^2 +9 s^3) z^2 {\varphi'}^2
} \\ \displaystyle{ \phantom{1234}
- (33 +18 s +3 s^2) s z \varphi \varphi'
-9 \varphi^2
-8
=0.
}
\end{array}
\right.
\label{eqoptim3}
\end{eqnarray}

An equivalent, shorter expression is obtained by suppressing
the term $\varphi^2$ \cite[p.~46]{F2004},
\begin{eqnarray}
& & \left\lbrace
\begin{array}{ll}
\displaystyle{
z=x y^{-\mu},\ f=\left(\frac{x y}{z}\right)^{3/2} \varphi(z),
\hbox{ or }
z=y x^{-\mu},\ f=\left(\frac{x y}{z}\right)^{3/2} \varphi(z),
}
\\
\displaystyle{
\left[
16 \mu^2 (\mu-1) z^2 \varphi''
-8 \mu (3 \mu+1)(\mu+1) z \varphi'
+3 (3 \mu+1) (3 \mu-1) (\mu+1) \varphi
\right] \varphi'''
} \\ \displaystyle{ \phantom{1234}
-8 \mu (\mu-3) z {\varphi''}^2
+(\mu-3)(\mu+3) (\mu+1) \varphi' \varphi''
-8=0.
}
\end{array}
\right.
\label{eqoptim3bis}
\end{eqnarray}
As it results from the scaling solution (\ref{eqScalingSolution}),
the ODE (\ref{eqoptim3bis}) admits the particular zero-parameter solution
\begin{eqnarray}
& &
\forall \mu\ :\
\varphi=2 i \frac{\sqrt{2}}{3} z^{3/2},\
      f=2 i \frac{\sqrt{2}}{3} (x y)^{3/2}.
\label{eqoptim3bisPS}
\label{eqoptim3mu1PS}
\end{eqnarray}
For generic values of $(a,b)$,
this ODE is unfortunately outside the class
\begin{eqnarray}
& &
\varphi''' = \sum_{j=0}^{3} A_j(z,\varphi,\varphi') {\varphi''}^j,
\end{eqnarray}
an equation which for some $A_j$ can be linearized by a contact
transformation. 
However, there exist particular values of $\mu$ for which
the integration can be performed at least partially.
The invariance of (\ref{eqoptim3}) under $(a,b) \to (b,a)$
induces an invariance of (\ref{eqoptim3bis}) under $\mu \to \mu^{-1}$.

\begin{enumerate}
\item
For $\mu=0,1,-1,-2,-1/2$, a first integral $K$ is known,
\begin{eqnarray}
& & \left\lbrace \matrix{
\displaystyle{
\mu=0,\
K=-8 z -3 \varphi \varphi'' - 3 {\varphi'}^2
=\left[-\frac{4}{3} z^3 - \frac{3}{2} \varphi^2\right]'',
}
\hfill \cr
\displaystyle{
\mu=1,\
K=\hbox{any rational function of } a,b,c,
\hbox{ see } (\ref{eqoptim3mu1FirstIntegrals}),
}
\hfill \cr
\displaystyle{
\mu=-1,\
K=z + 2 z^2 {\varphi''}^2,
}
\hfill \cr
\displaystyle{
\mu=-2,\
K=-8 z
- {\varphi'}^2
-105 \varphi \varphi''
+ 112 z \varphi' \varphi''
-96 z^2 {\varphi''}^2,
}
\hfill \cr
\displaystyle{
\mu=-\frac{1}{2},\
K=-8 z^2+\frac{15}{4}
  \left(z \varphi \varphi'' + z {\varphi'}^2 - \varphi \varphi'\right)
-10 z^2 \varphi' \varphi'' -6 z^3 {\varphi''}^2.
}
\hfill \cr
} \right.
\end{eqnarray}

\item
For $\mu=1$, the third order equation \cite[Eq.~(31) p.~46]{F2004},
\begin{eqnarray}
& & \left\lbrace \matrix{
\displaystyle{
z=x/y,\ f=y^3 \varphi(z),\
\hbox{ or }
z=y/x,\ f=x^3 \varphi(z),\
}
\hfill \cr
\displaystyle{
 2(3 \varphi - 2 z \varphi') \varphi''' + 2 z {\varphi''}^2
  - 2 \varphi' \varphi'' -1=0,
} \hfill
\cr} \right.
\label{eqoptim3mu1}
\end{eqnarray}
is linearizable since its derivative factorizes into
\begin{eqnarray}
& &
2 (3 \varphi - 2 z \varphi') \varphi''''=0,
\end{eqnarray}
so its general solution is
\begin{eqnarray}
& &
\varphi=\alpha z^3 + 3 \beta z^2 + 3 \gamma z + \delta,\
36 (\alpha \delta - \beta \gamma) -1=0,
(\alpha,\beta,\gamma) \hbox{ arbitrary}.
\label{eqoptim3mu1GS}
\end{eqnarray}
It is interesting to notice that, knowing the three first integrals $a,b,c$,
\begin{eqnarray}
& & \left\lbrace \matrix{
\displaystyle{
12 a=\frac{1 + 2 \varphi' \varphi'' - 2 z {\varphi''}^2}
           {3 \varphi -2 z \varphi'}
    = 2 \varphi''',
}
\hfill \cr
\displaystyle{
4 b=\frac{-z + 6 \varphi \varphi'' - 6 z \varphi' \varphi''
            +2 z^2 {\varphi''}^2}
           {3 \varphi -2 z \varphi'}
   =2 \varphi'' - 2 z \varphi''',
}
\hfill \cr
\displaystyle{
4 c=\frac{z^2 + 12 \varphi \varphi' - 8 z {\varphi'}^2
            -12 z \varphi \varphi'' + 10 z^2 \varphi' \varphi''
            -2 z^3 {\varphi''}^2}
           {3 \varphi -2 z \varphi'}
   =4 \varphi' - 4 z \varphi'' + 2 z^2 \varphi''',
}
\hfill \cr
} \right.
\label{eqoptim3mu1FirstIntegrals}
\end{eqnarray}
there exists no first integral
which would be polynomial in $(\varphi,\varphi',\varphi'')$.

\item
For $\mu=-1$, the ODE reduces to a linear equation for ${\varphi''}^2$,
identical to the particular case $r_1=r_2=s_1=s_2=0$ of
the reduction (\ref{eqoptim4}) given below.

\item
For $\mu=3,\ 1/3$ and $\mu=-3,\ -1/3$ respectively,
the ODE is just the subcase $a=0$ of
the reductions (\ref{eqoptim5or6order3}) and (\ref{eqoptim7or8}) given below.

\item
For $\mu=2,1/2$,
two rational solutions for $\varphi^2$ can be obtained,
\begin{eqnarray}
& & \left\lbrace \matrix{
\displaystyle{
\mu=2,\ \varphi=\frac{2}{15 c} (z-c)^{5/2},\
f=\frac{2 y^2}{15 c} \left(\frac{x}{y}-c y \right)^{5/2},
}
\hfill \cr
\displaystyle{
\mu=\frac{1}{2},\ \varphi=\frac{2}{15 c} z^{-1/2} (1-c z^2)^{5/2},\
f=\frac{2 x^2}{15 c} \left(\frac{y}{x}-c x \right)^{5/2},
}
\hfill \cr
} \right.
\label{eqoptim3subcase2and1over2old}
\end{eqnarray}
and
\begin{eqnarray}
& & \left\lbrace \matrix{
\displaystyle{
\mu=2,\ \varphi=2 i \frac{\sqrt{2}}{3} z^{3/2} (1-c z),\
f=\frac{2 i \sqrt{2}}{3} (x y)^{3/2}\left(1-\frac{c x}{y^2}\right)
}
\hfill \cr
\displaystyle{
\mu=\frac{1}{2},\ \varphi=2 i \frac{\sqrt{2}}{3} z^{-1/2} (z^2-c),\
f=\frac{2 i \sqrt{2}}{3} (x y)^{3/2}\left(1-\frac{c y}{x^2}\right)
}
\hfill \cr
} \right.
\label{eqoptim3subcase2and1over2new}
\end{eqnarray}
in which $c$ is arbitrary.

The first solution (\ref{eqoptim3subcase2and1over2old})
represents the octahedron solution $B_3$ of Dubrovin,
see \cite[p.~41]{F2004}.

The second solution $f$ extrapolates the scaling solution
(\ref{eqScalingSolution}).

\item
For $\mu=5/3,\ 3/5$,
one rational solution exists,
which depends on one arbitrary parameter $c$,
\begin{eqnarray}
& & \left\lbrace \matrix{
\displaystyle{
\mu=\frac{5}{3},\ \varphi=\frac{c}{6} z^3 + \frac{1}{24 c},\
}
\hfill \cr
\displaystyle{
\mu=\frac{3}{5},\ \varphi=\frac{c}{6 z} + \frac{z^4}{24 c}.
}
\hfill \cr
} \right.
\label{eqoptim3subcase5over3}
\end{eqnarray}
This represents the tetrahedron solution $A_3$ of Dubrovin,
see \cite[p.~41]{F2004}.

\end{enumerate}

\subsection{Reduction with the generator $\bfv_3-\bfv_4+a\bfv_5+ b \bfv_{10}$}

This reduction to a nonautonomous ODE,
\begin{eqnarray}
& & \left. \matrix{
\displaystyle{
z=xy,\ f=\varphi(z)+(a z + b) \log x,
}
\hfill \cr
}
\right.
\label{eqoptim4bad}
\end{eqnarray}
can be defined more symmetrically as \cite[p.~45, Example 2]{F2004}
\begin{eqnarray}
& & \left\lbrace
\begin{array}{ll}
\displaystyle{
z=xy,\
f=\varphi(z) + (r_1 z + r_2) \log x + (s_1 z + s_2) \log y,
}
\\
\displaystyle{
z^2 {\varphi''}^2 +(r_1+s_1) z \varphi'' -(r_2+s_2) \varphi''
- \frac{r_1 s_2 + r_2 s_1}{z} + \frac{r_2 s_2}{z^2} + \frac{z}{2}
+ \frac{(r_1+s_1)^2}{4} + k_0 =0,
}
\end{array}
\right.
\label{eqoptim4}
\end{eqnarray}
in which $k_0$ is a constant of integration.
Its general solution is obtained by quadratures,
\begin{eqnarray}
& & \left\lbrace
\begin{array}{ll}
\displaystyle{
\varphi=k_1 z + k_2 - \frac{r_1+s_1}{2} (z \log z -z)
- \frac{r_2+s_2}{2} \log z
}
\\
\displaystyle{
\phantom{1234}
\pm \int \D z \int \D z \frac{\sqrt{
-2 z^3 - 4 k_0 z^2 - 2(r_1-s_1)(r_2-s_2) z+(r_2-s_2)^2}}{2 z^2},
}
\\
\displaystyle{
f=-s_2 \log x - r_2 \log y +\frac{r_1-s_1}{2} x y \log\frac{x}{y}
\pm \int \D z \int \D z \frac{\sqrt{\dots}}{2 z^2},
}
\end{array}
\right.
\label{eqoptim4GS}
\end{eqnarray}
and it generically involves elliptic integrals.
A particular solution is
\begin{eqnarray}
& & \left\lbrace
\begin{array}{ll}
\displaystyle{
f=2 i \frac{\sqrt{2}}{3} (x y)^{3/2} + c x y \log \frac{x}{y},
}
\\
\displaystyle{
F= i \frac{\sqrt{2}}{8} x^{-1/2} y^{3/2}
\left(4 c^2 \log\frac{x}{y} - x y \right)
+\frac{c y^2}{4} + \frac{c^3 y}{x}
+\left(\frac{c^3 y}{x}-\frac{c y^2}{2}\right) \log\frac{x}{y},
}
\end{array}
\right.
\label{eqoptim4PS}
\end{eqnarray}
which is another extrapolation of
the scaling solution (\ref{eqScalingSolution}).

\subsection{Reduction with the generator
 $3 \bfv_3 + \bfv_4 + a \bfv_7$ or $\bfv_3 + 3 \bfv_4 + a \bfv_6$}

These two generators define a reduction to the same
nonautonomous ODE,
\begin{eqnarray}
& & \left\lbrace
\begin{array}{ll}
\displaystyle{
z=x y^{-3},\ f= y^6 \varphi(z) - \frac{a}{4} y^2,
\hbox{ or }
z=y x^{-3},\ f= x^6 \varphi(z) - \frac{a}{4} x^2,
}
\\
\displaystyle{
12(3 z^2 \varphi'' -8 z \varphi' +10 \varphi) \varphi'''-1=0,
}
\end{array}
\right.
\label{eqoptim5or6order3}
\end{eqnarray}
which a linear transformation can make second order in $\varphi'$,
\begin{eqnarray}
& & \left\lbrace
\begin{array}{ll}
\displaystyle{
z=x y^{-3},\ f= x^2 \varphi(z) - \frac{a}{4} y^2,
\hbox{ or }
z=y x^{-3},\ f= y^2 \varphi(z) - \frac{a}{4} x^2,
}
\\
\displaystyle{
\left[36 z^6 \varphi'' + 48 z^5 \varphi' \right] \varphi'''
+ 216 z^5 {\varphi''}^2
+504 z^4 \varphi' \varphi''
+288 z^3 {\varphi'}^2
-1
=0.
}
\end{array}
\right.
\label{eqoptim5or6order2}
\end{eqnarray}
Since $f$ is defined up to an arbitrary additive second degree polynomial,
the reduced ODE does not depend on $a$,
and this case is identical to the case $\mu=3,1/3$ of (\ref{eqoptim3bis}),
in which no solution is known other than (\ref{eqoptim3bisPS}).

\subsection{Reduction with the generator
 $-3 \bfv_3 +   \bfv_4 + a \bfv_8$ or $\bfv_3 - 3 \bfv_4 + a \bfv_9$}

These two generators define a reduction to the same
second order, nonautonomous ODE for $\varphi'$,
\begin{eqnarray}
& & \left\lbrace
\begin{array}{ll}
\displaystyle{
z=x y^3,\ f= x \varphi(z) - \frac{a}{3} x \log x,
\hbox{ or }
z=y x^3,\ f= y \varphi(z) - \frac{a}{3} y \log y,\
a \not=0,
}
\\
\displaystyle{
\left[-72 z^4 \varphi'' -84 z^3 \varphi' + 9 a z^2 \right] \varphi'''
- 234 z^3 {\varphi''}^2
} \\ \displaystyle{ \phantom{1234}
-324 z^2 \varphi' \varphi''
+ 18 a z  \varphi''
- 72 z {\varphi'}^2
+ 2 a \varphi'
-1=0,
}
\end{array}
\right.
\label{eqoptim7or8}
\end{eqnarray}
but, with $a \not=0$, we could not find any solution to this ODE.

\subsection{Reduction with the generator
 $a \bfv_2 + b \bfv_3$ or $a \bfv_1 + b \bfv_4$}

They lead to the same autonomous ODE,
\begin{eqnarray}
& & \left\lbrace
\begin{array}{ll}
\displaystyle{
z=b x - a \log y,\ f=y^{3/2} \varphi(z),
\hbox{ or }
z=a y - b \log x,\ f=x^{3/2} \varphi(z),\
a b \not=0,
}
\\
\displaystyle{
\left[16 a^2 \varphi'' - 8 a \varphi' +3 \varphi \right] \varphi'''
- 24 a {\varphi''}^2 + 9 \varphi' \varphi'' + 8 b^{-3}=0.
}
\end{array}
\right.
\label{eqoptim9or10}
\end{eqnarray}
We could not find a particular solution for this ODE.


\subsection{Reduction with the generator $-a \bfv_1 + b \bfv_2                       + c \bfv_5 + d \bfv_6 + e \bfv_7$}

The reduced ODE is autonomous and linear \cite[p.~44, Example 1]{F2004},
\begin{eqnarray}
& & \left\lbrace
\begin{array}{ll}
\displaystyle{
 z=b x + a y,\ f= \varphi(z) + c_3 x^3 + c_2 x^2 y + c_1 x y^2 + c_0 y^3,\
a b \not=0,
}
\\
\displaystyle{
c=-2 a c_2+2 b c_1,d=-3 a c_3+b c_2,e=-a c_1+3 b c_0,
}
\\
\displaystyle{
2 \left(3 a^3 c_3 - a^2 b c_2 - a b^2 c_1 + 3 b^3 c_0\right) \varphi'''
+ 36 c_0 c_3 - 4 c_1 c_2 -1=0.
}
\end{array}
\right.
\label{eqoptim11}
\end{eqnarray}
and the solution $f(x,y)$ (always defined up to an arbitrary polynomial of
degree two in $(x,y)$)
is identical to that defined by Eq.~(\ref{eqoptim3mu1GS}),
i.e.~the third degree polynomial depending on three
arbitrary independent constants,
\begin{eqnarray}
& &
f(x,y)=\alpha x^3 + 3 \beta x^2 y + 3 \gamma x y^2 + \delta y^3,\
36 (\alpha \delta - \beta \gamma) -1=0.
\label{eqoptim11GS}
\end{eqnarray}

\section {Summary of solutions}
\label{sectionSummary}

The explicit solutions to (\ref{eqPDEFerapontov})
are summarized in Table \ref{Table0}.
This table does not include the reductions for which no solution could
be found.
The too long expression for the ``icosa$'$'' solution is,
\begin{eqnarray}
& & \left\lbrace
\begin{array}{ll}
\displaystyle{
F(y,t)=\frac{   \irrec^2 x^6 T}{4}
     + \frac{29 \irrec^3 x^5 T^4}{24}
      +\frac{29 \irrec^4 x^4 T^7}{30}
      +\frac{   \irrec^5 x^3 T^{10}}{10}
      +\frac{ 3 \irrec^6 x^2 T^{13}}{80}
      +\frac{   \irrec^8     T^{19}}{3040},
}
\\
\displaystyle{
f(x,y)=
  \frac{4 \irrec^4 x^2 T^9}{45}
+ \frac{7 \irrec^3 x^3 T^6}{30}
+ \frac{  \irrec^2 x^4 T^3}{6}
+ \frac{  \irrec   x^5    }{60},
}
\\
\displaystyle{
y=\irrec x^2 T + \frac{\irrec^2 x T^4}{2},\
t=\frac{\irrec x^3}{3} + \irrec^2 x^2 T^3 +\frac{\irrec^4 T^9}{36}.
}
\end{array}
\right.
\label{eqicosaprime}
\end{eqnarray}



\begin{table}[h] 
\caption[garbage]{
Summary of solutions $F(y,t),f(x,y)$ of the equations
(\ref{eqPDEDubrovin}), (\ref{eqPDEFerapontov}).
A3, B3, H3 label solutions linked to regular polyhedra \cite{DubLNM},
Dubn solutions found by Dubrovin \cite{DubLNM},
Fn additional solutions listed in \cite[p.~41]{F2004},
and Nn solutions apparently new.
A prime (') labels the solution deduced by permuting $x$ and $y$ in $f$.
A blank field in column ``Eq'' indicates a solution not arising from
a known reduction.
The irrelevant constant $\irrec$ reflects the scaling invariance
and can be set to $1$.
$P_n$ denotes a polynomial of degree $n$.
}
\vspace{0.2truecm}
\begin{center}
\begin{tabular}{| l | l | l | l | l |}
\hline 
Label
&
$F(y,t)$ 
&
$f(x,y)$ 
&
Eq
&
Link $(t,x,y)$
\\ \hline   \hline 
&
$\displaystyle{\frac{y^4}{8 t}}$
&
$\displaystyle{2 i \frac{\sqrt{2}}{3}(x y)^{3/2} }$
&
(\ref{eqScalingSolution})
&
$\displaystyle{x=-\frac{y^3}{2 t^2}}$
\\ \hline 
F1
&
$\displaystyle{i \sqrt{2} y^{5/2} x^{-15/2} \lambda^{-5} P_{8}(x)}$
&
$\displaystyle{2 i \frac{\sqrt{2}}{3}(x y)^{3/2} \lambda,\
\lambda^2=1+\frac{\alpha}{x}+\frac{\beta}{x^2}+\frac{\gamma}{x^3}}$
&
(\ref{eqoptim1or2})
&
$\displaystyle{t=i \sqrt{2} y^{3/2} x^{-9/2} \lambda^{-3} P_{4}(x)}$
\\ \hline 
F1'
&
$\displaystyle{
\frac{4 \alpha \gamma-\beta^2 +12 \gamma y +6 \beta y^2 +4 \alpha y^3 +3 y^4}
{24 t}
}$
&
$\displaystyle{2 i \frac{\sqrt{2}}{3}(x y)^{3/2} \lambda,\
\lambda^2=1+\frac{\alpha}{y}+\frac{\beta}{y^2}+\frac{\gamma}{y^3}}$
&
(\ref{eqoptim1or2})
&
$\displaystyle{x=-\frac{\lambda^2 y^3}{t^2}}$
\\ \hline   \hline 
F2
&
$\displaystyle{}$ 
&
$\displaystyle{
{\displaystyle{
\frac{r_1-s_1}{2} x y \log\frac{x}{y}
-s_2 \log x - r_2 \log y
}}
\atop
{\displaystyle{\pm \int \D z \int \D z \frac{\sqrt{\dots}}{2 z^2},\ z=xy}}
}$
&
(\ref{eqoptim4GS})
&
$t=f_{xx}$
\\ \hline  
F3
&
(\ref{eqoptim4PS})
&
$\displaystyle{2 i \frac{\sqrt{2}}{3}(x y)^{3/2} + c x y \log\frac{x}{y}}$
&
(\ref{eqoptim4PS})
&
$\displaystyle{
t=i \frac{\sqrt{2}}{2} \left(\frac{y^3}{x}\right)^{1/2} + c \frac{y}{x}}$
\\ \hline   \hline 
F4
&
$\displaystyle{
\frac{12 (\beta^2-\alpha \gamma) y^3 -6 \beta y^2 t + y t^2 + 2 \gamma t^3}
     {12 \alpha}
}$
&
$\displaystyle{
{\alpha x^3 + 3 \beta x^2 y + 3 \gamma x y^2 + \delta y^3}
\atop
{36 (\alpha \delta - \beta \gamma) -1=0}
}$
&
$\displaystyle{{(\ref{eqoptim3mu1GS})} \atop {(\ref{eqoptim11GS})}}$
&
$t=6 (\alpha x + \beta y)$
\\ \hline   \hline 
$\displaystyle{\hbox{octa}\atop\hbox{B3}}$
&
$\displaystyle{
 \frac{\irrec y^3 t}{3} + \frac{2 \irrec^2 y^2 t^3}{3}+\frac{8 \irrec^4 t^7}{105}}$
&
$\displaystyle{\frac{2 y^2}{15 \irrec} \left(\frac{x}{y}-\irrec y \right)^{5/2}
}$
&
(\ref{eqoptim3subcase2and1over2old})
&
$\displaystyle{x=4 \irrec^2 y t^2 + \irrec y^2}$
\\ \hline 
octa'
&
$\displaystyle{
\displaystyle{
 \frac{\lambda^{11}}{528 \irrec^3}
+\frac{\lambda^{7} x^2}{6 \irrec}
+\frac{\lambda^{5} x^3}{3}
+\frac{7 \irrec \lambda^{3} x^4}{3}
}
\atop
\displaystyle{
+\frac{4 \irrec^2 \lambda x^5}{3}
}}$
&
$\displaystyle{
\frac{2 x^2}{15 \irrec} \left(\frac{y}{x}-\irrec x \right)^{5/2},\
y=\irrec x^2 + x \lambda^2
}$
&
(\ref{eqoptim3subcase2and1over2old})
&
$\displaystyle{t=2 \irrec \lambda x^2 + \frac{\lambda^5}{10 \irrec}}$
\\ \hline   \hline 
N1 
&
$\displaystyle{
\displaystyle{
\frac{i \sqrt{2}}{24}x^{1/2} y^{-7/2}\times
}
\atop
\displaystyle{
\left(
\frac{25 \irrec^3 x^3}{7}-5 \irrec^2 x^2 y^2 -7 \irrec x y^4 - 3 y^6\right)
}}$
&
$\displaystyle{
2i \frac{\sqrt{2}}{3} (x y)^{3/2}\left(1-\frac{\irrec x}{y^2}\right)
}$
&
(\ref{eqoptim3subcase2and1over2new})
&
$\displaystyle{t=\frac{i \sqrt{2}}{2 (x y)^{1/2}} (y^2-5 \irrec x)}$
\\ \hline  
N1' 
&
$\displaystyle{
\displaystyle{\frac{i \sqrt{2}}{24}x^{-11/2} y^{5/2}\times
}
\atop
\displaystyle{
\left(
\frac{125 \irrec^3 y^3}{11}-25 \irrec^2 x^2 y^2 +5 \irrec x^4 y -3 x^6\right)}
}$
&
$\displaystyle{
2i \frac{\sqrt{2}}{3} (x y)^{3/2}\left(1-\frac{\irrec y}{x^2}\right)
}$
&
(\ref{eqoptim3subcase2and1over2new})
&
$\displaystyle{t=\frac{i \sqrt{2}}{2} x^{-5/2} y^{3/2} (x^2-\irrec y)}$
\\ \hline  \hline 
$\displaystyle{\hbox{tetra}\atop\hbox{A3}}$
&
$\displaystyle{\frac{y^2 t^2}{4 \irrec} + \frac{t^5}{60 \irrec^2}}$
&
$\displaystyle{\frac{\irrec x^3}{6 y} + \frac{y^4}{24 \irrec}}
$
&
(\ref{eqoptim3subcase5over3})
&
$x=t y / \irrec$ 
\\ \hline 
tetra'
&
$\displaystyle{
\frac{x^3 y}{3 \irrec} + \frac{3 \irrec y^4}{8 x^2}+\frac{\irrec^3 y^7}{28 x^7}
}$
&
$\displaystyle{\frac{\irrec y^3}{6 x} + \frac{x^4}{24 \irrec}}
$
&
(\ref{eqoptim3subcase5over3})
&
$\displaystyle{t=\frac{\irrec y^3}{3 x^3} + \frac{x^2}{2 \irrec}}$
\\ \hline \hline 
Dub1
&
$\displaystyle{
 \frac{e^{2 \irrec t}}{8 \irrec^3}
+\frac{y^2 e^{\irrec t}}{2 \irrec}
-\frac{\irrec y^4}{48}}$
&
$\displaystyle{
\frac{\irrec x y^3}{12} - \frac{x^2}{2 \irrec} \log \frac{x}{y}
   -\frac{3 x^2}{4 \irrec}
}$
&
&
$\displaystyle{x=y e^{\irrec t}}$
\\ \hline 
Dub1'
&
$\displaystyle{
{\displaystyle{
 \frac{y^4}{32 \irrec^3 x^4}(4 \log\frac{y}{x}-3)
}}
\atop
{\displaystyle{
+\left(\frac{y^3}{ 8 \irrec x}+\frac{\irrec x^2 y^2}{16}\right)
  (2 \log\frac{y}{x}+3)
}}
}$
&
$\displaystyle{
\frac{\irrec y x^3}{12} - \frac{y^2}{2 \irrec} \log \frac{y}{x}
   -\frac{3 y^2}{4 \irrec}
}$
&
&
$\displaystyle{t=\frac{\irrec x y}{2} + \frac{y^2}{2 \irrec x^2}}$
\\ \hline  \hline 
Dub2
&
$\displaystyle{-\frac{\irrec y^4}{24} + \frac{y}{\irrec} e^{\irrec t}}$
&
$\displaystyle{
\frac{\irrec x y^3}{6} + \frac{x^2}{2 \irrec} \log x -\frac{3 x^2}{4 \irrec}
}$
&
&
$\displaystyle{x=e^{\irrec t}}$
\\ \hline 
Dub2'
&
$\displaystyle{\frac{t^2 \log y}{2 \irrec}}$
&
$\displaystyle{
\frac{\irrec x y^3}{6} + \frac{x^2}{2 \irrec} \log x -\frac{3 x^2}{4 \irrec}
}$
&
&
$\displaystyle{x=\frac{t}{\irrec y}}$
\\ \hline  \hline 
$\displaystyle{\hbox{icosa}\atop\hbox{H3}}$
&
$\displaystyle{
\frac{\irrec y^3 t^2}{6} + \frac{\irrec^2 y^2 t^5}{20}
+ \frac{\irrec^4 t^{11}}{3960}
}$
&
$\displaystyle{\frac{\irrec^2 y^4 t^3}{6} + \frac{7 \irrec^3 y^3 t^6}{30}
              + \frac{4 \irrec^4 y^2 t^9}{45} + \frac{\irrec y^5}{60}
}$
&
&
$\displaystyle{x=\irrec y^2 t + \frac{\irrec^2 y t^4}{2}}$
\\ \hline 
icosa'
&
(\ref{eqicosaprime})
&
$\displaystyle{
  \frac{4 \irrec^4 x^2 T^9}{45}
+ \frac{7 \irrec^3 x^3 T^6}{30}
+ \frac{  \irrec^2 x^4 T^3}{6}
+ \frac{  \irrec   x^5    }{60}
}$
&
&
$\displaystyle{
\displaystyle{y=\irrec x^2 T + \frac{\irrec^2 x T^4}{2}}
\atop
\displaystyle{
t=\frac{\irrec x^3}{3} + \irrec^2 x^2 T^3 +\frac{\irrec^4 T^9}{36}}
}$
\\ \hline 
\end{tabular}
\end{center}
\label{Table0}
\end{table}

\section{Conclusion}

Finding additional solutions to the obtained reductions
could generate algebraic solutions of the sixth Painlev\'e equation $\PVI$
\cite{DubLNM},
in which the four monodromy exponents of $\PVI$
could depend on one arbitrary constant,
like in some particular cases (tetrahedron and octahedron solutions)
found by Kitaev \cite{KitaevP6cube}.
In particular, the two solutions labeled N1 and N1' in Table \ref{Table0}
obey the quasi-homogeneity condition (WDVV3) recalled in the introduction.
This question is currently under investigation.

\section*{Acknowledgments}

We warmly thank Evgueni Ferapontov and Zhang You-jin for 
enlightening discussions,
and one referee for suggestions to greatly improve the introduction.


\end{document}